\documentclass[10pt]{article}

\usepackage{amscd,amsmath, amssymb, fancyhdr, mathbbol}

\numberwithin{equation}{section}

\newcommand{\version}{version 3.0,\ \ March 2, 2017}

\def\eqref#1{(\ref{#1})}

\newcommand{\arrow}{{\:\longrightarrow\:}}
\newcommand{\Z}{{\Bbb Z}}
\def\C{{\Bbb C}}

\def\1{\sqrt{-1}\:}

\newcommand{\cntrct}                
{\hspace{2pt}\raisebox{1pt}{\text{$\lrcorner$}}\hspace{2pt}}

\renewcommand{\phi}{\varphi}
\renewcommand{\epsilon}{\varepsilon}
\renewcommand{\geq}{\geqslant}
\renewcommand{\leq}{\leqslant}

\newcommand{\Aut}{\operatorname{Aut}}

\newcounter{Mycounter}[section]
\newcounter{lemma}[section]
\renewcommand{\thelemma}{{Lemma \thesection.\arabic{lemma}}}
\newcommand{\lemma}{%
    \setcounter{lemma}{\value{Mycounter}}
    \refstepcounter{lemma}
    \stepcounter{Mycounter}
    {\noindent \bf \thelemma:\ }}

\newcounter{claim}[section]
\newcounter{sublemma}[section]
\newcounter{corollary}[section]
\newcounter{theorem}[section]
\renewcommand{\thetheorem}{{Theorem \thesection.\arabic{theorem}}}
\newcommand{\theorem}{%
    \setcounter{theorem}{\value{Mycounter}}
    \refstepcounter{theorem}
    \stepcounter{Mycounter}
    {\noindent \bf \thetheorem:\ }}

\newcounter{conjecture}[section]
\renewcommand{\theconjecture}{{Conjecture \thesection.\arabic{conjecture}}}
\newcommand{\conjecture}{%
    \setcounter{conjecture}{\value{Mycounter}}
    \refstepcounter{conjecture}
    \stepcounter{Mycounter}
    {\noindent \bf \theconjecture:\ }}

\newcounter{proposition}[section]
\renewcommand{\theproposition}
      {{Proposition \thesection.\arabic{proposition}}}
\newcommand{\proposition}{%
    \setcounter{proposition}{\value{Mycounter}}
    \refstepcounter{proposition}
    \stepcounter{Mycounter}
    {\noindent \bf \theproposition:\ }}

\newcounter{definition}[section]
\renewcommand{\thedefinition}
      {{Definition~\thesection.\arabic{definition}}}
\newcommand{\definition}{%
    \setcounter{definition}{\value{Mycounter}}
    \refstepcounter{definition}
    \stepcounter{Mycounter}
    {\noindent \bf \thedefinition:\ }}

\newcounter{example}[section]
\newcounter{remark}[section]
\renewcommand{\theremark}{{Remark \thesection.\arabic{remark}}}
\newcommand{\remark}{%
    \setcounter{remark}{\value{Mycounter}}
    \refstepcounter{remark}
    \stepcounter{Mycounter}
    {\noindent \bf \theremark:\ }}

\newcounter{problem}[section]
\newcounter{question}[section]
\makeatletter

\@addtoreset{equation}{section} \@addtoreset{footnote}{section}
\makeatother

\def\blacksquare{\hbox{\vrule width 5pt height 5pt depth 0pt}}
\def\endproof{\blacksquare}

\begin{document}
\begin{center}
{\LARGE\bf
Algebraic non-hyperbolicity of hyperk\"ahler manifolds 
with Picard rank greater than one\\[4mm]
}

Ljudmila Kamenova, Misha Verbitsky\footnote{Partially supported 
by RSCF grant 14-21-00053 within AG Laboratory NRU-HSE.}

\end{center}

{\small \hspace{0.10\linewidth}
\begin{minipage}[t]{0.85\linewidth}
{\bf Abstract.} 
A projective manifold is algebraically hyperbolic if the degree of any curve
is bounded from above by its genus times a constant, which is independent 
from the curve.
This is a property which follows from Kobayashi hyperbolicity.
We prove that hyperk\"ahler manifolds are not algebraically 
hyperbolic when the Picard rank is at least 3, or if
the Picard rank is 2 and the SYZ conjecture on existence
of Lagrangian fibrations is true.
We also prove that if the automorphism group of a hyperk\"ahler 
manifold is infinite then it is algebraically non-hyperbolic. 
\end{minipage}
}

\section{Introduction} 

In \cite{_Verbitsky:ergodic_} M. Verbitsky proved that all hyperk\"ahler 
manifolds are Kobayashi non-hyperbolic. It is interesting to inquire if 
projective hyperk\"ahler manifolds are also algebraically 
non-hyperbolic (\ref{_alge_hype_Definition_}). 
For a given projective manifold algebraic non-hyperbolicity 
implies Kobayashi non-hyperbolicity. We prove algebraic non-hyperbolicity 
for projective hyperk\"ahler manifolds with infinite group of 
automorphisms. 

\hfill

\theorem
Let $M$ be a projective hyperk\"ahler manifold with infinite automorphism 
group. Then $M$ is algebraically non-hyperbolic. 

\hfill

If a projective hyperk\"ahler manifold has Picard rank at least three, 
we show that it is algebraically non-hyperbolic. For the case when the 
Picard rank equals to two we need an extra assumption in order to prove 
algebraic non-hyperbolicity. The SYZ conjecture states that a nef 
parabolic line bundle on a hyperk\"ahler manifold gives rise to a 
Lagrangian fibration (\ref{syz}). 

\hfill

\theorem
Let $M$ be a projective hyperk\"ahler manifold with Picard rank $\rho$. 
Assume that either $\rho > 2$, or $\rho =2$ and 
the SYZ conjecture holds. 
Then $M$ is algebraically non-hyperbolic.

\section{Basic notions}

\definition
A {\em hyperk\"ahler manifold of maximal holonomy} 
(or {\em irreducible holomorphic symplectic}) 
manifold 
$M$ is a compact complex K\"ahler manifold with $\pi_1(M)=0$ and 
$H^{2,0}(M)=\C \sigma$, where $\sigma$ is everywhere non-degenerate. 
From now on we would tacitly assume that hyperk\"ahler manifolds are
of maximal holonomy.

\hfill

Due to results of Matsushita, holomorphic maps from hyperk\"ahler 
manifolds are quite restricted. 

\hfill

\begin{theorem} (Matsushita, \cite{_Matsushita:fibred_}) \label{fibr}
Let $M$ be a hyperk\"ahler manifold and $f\colon M\rightarrow B$ a 
proper surjective morphism with a smooth base $B$. Assume that $f$ has 
connected fibers and $0 < \dim B < \dim M$. Then $f$ is Lagrangian 
and $\dim_{\C} B = n$, where $\dim_{\C} M = 2n$. 
\end{theorem}

\hfill

Following \ref{fibr}, we call the surjective morphism 
$f\colon M\rightarrow B$ a {\em Lagrangian fibration} on the 
hyperk\"ahler manifold $M$. 
A dominant map $f\colon M\dashrightarrow B$ 
is a {\em rational Lagrangian fibration} if there exists a 
birational map $\phi \colon M \dasharrow M'$ between hyperk\"ahler 
manifolds such that the composition 
$f \circ \phi^{-1} \colon M' \rightarrow B$ is 
a Lagrangian fibration. J.-M. Hwang proved that if the base $B$ of a 
hyperk\"ahler Lagrangian fibration is smooth, then $B \cong {\mathbb P}^n$ 
(see \cite{Hwang}). 

\hfill

\definition
Given a hyperk\"ahler manifold $M$, there is a non-degenerate primitive 
form $q$ on $H^2(M,\Z)$, called the {\em Beauville-Bogomolov-Fujiki form} 
(or the {\em ``BBF form''} for short), 
of signature $(3,b_2-3)$, and satisfying the {\em Fujiki relation} 
$$\int_M \alpha^{2n}=c\cdot q(\alpha)^n\qquad\text{for }\alpha \in 
H^2(M,\Z),$$ 
with $c>0$ a constant depending on the topological type of $M$. 
This form generalizes the intersection pairing on K3 surfaces. 
A detailed description of the form can be found in 
\cite{_Beauville_}, \cite{_Bogomolov:defo_} and \cite{_Fujiki:HK_}. 

\hfill

Notice that given a Lagrangian fibration $f\colon M \rightarrow {\mathbb P}^n$, 
if $h$ is the hyperplane class on $\mathbb P^n$, and $\alpha=f^*h$, 
then $\alpha$ belongs to the birational K\"ahler cone of $M$ and 
$q(\alpha)=0$. The following SYZ conjecture states that the converse 
is also true. 

\hfill

\conjecture (Tyurin, Bogomolov, Hassett-Tschinkel,
Huybrechts, Sawon) \label{syz}
If $L$ is a line bundle on a hyperk\"ahler manifold $M$ with $q(L)=0$, 
and such that $c_1(L)$ belongs to the birational K\"ahler cone of $M$, 
then $L$ defines a rational Lagrangian fibration.

\hfill

For more reference on this conjecture, please see
\cite{_Hassett_Tschinkel:SYZ_conj_}, 
\cite{_Sawon_}, \cite{_Huybrechts:lec_} and \cite{_Verbitsky:SYZ_}.
This conjecture is known for deformations of Hilbert schemes of points on 
K3 surfaces 
(Bayer--Macr\`i \cite{_Bayer_Macri_}; Markman \cite{_Markman:SYZ_}), 
and for deformations of the generalized Kummer varieties $K_n(A)$ 
(Yoshioka \cite{_Yoshioka_}).







\hfill

In \cite{_Verbitsky:ergodic_} M. Verbitsky proved that all hyperk\"ahler 
manifolds are Kobayashi non-hyperbolic. In \cite{klv} together with 
S. Lu we proved that the Kobayashi pseudometric vanishes 
identically for K3 surfaces and for hyperk\"ahler manifolds deformation 
equivalent to Lagrangian fibrations under some mild assumptions. 
In \cite{Demailly} Demailly introduced the following notion. 

\hfill

\definition\label{_alge_hype_Definition_}
A projective manifold $M$ is {\em algebraically hyperbolic} 
if for any Hermitian metric $h$ on $M$ there exists a constant $A>0$ 
such that for any holomorphic map $\phi \colon C \rightarrow M$ 
from a curve of genus $g$ to $M$ we have that 
$2g-2 \geq A \int_C \phi^\ast \omega_h,$ where $\omega_h$ is the K\"ahler 
form of $h$. 

\hfill

In this paper all varieties we 
consider are smooth and projective. For projective varieties, 
Kobayashi hyperbolicity implies algebraic hyperbolicity (\cite{Demailly}). 
Here we explore non-hyperbolic 
properties of projective hyperk\"ahler manifolds. Algebraic non-hyperbolicity
implies Kobayashi non-hyperbolicity. 



\section{Main Results}

\proposition\label{_hk_with_Lagra_non_h_Proposition_}
Let $M$ be a hyperk\"ahler manifold admitting a (rational) Lagrangian 
fibration. Then $M$ is algebraically non-hyperbolic. 

\hfill

\begin{proof}
We use the fact that the fibers of a Lagrangian fibrations are 
abelian varieties (\cite{_Matsushita:fibred_}). The isogeny self-maps 
on an abelian variety provide curves of fixed genus and arbirary large 
degrees, and therefore they are algebraically non-hyperbolic. 

\hfill

An alternative way of proving this proposition is by using the 
following result whose proof was suggested by Prof. K. Oguiso. 

\hfill

\lemma \label{ratlc}
If a hyperk\"ahler manifold $M$ admits a Lagrangian fibration, 
then there exists a rational curve on $M$. 

\hfill

Indeed, 
in \cite{HO} J.-M. Hwang and K. Oguiso give a Kodaira-type classification 
of the general singular fibers of a holomorphic Lagrangian fibration. 
All of the general singular fibers are covered by rational curves. 
The locus of singular fibers is non-empty (e.g., Proposition 4.1 in 
\cite{Hwang}), and therefore there is a rational curve on $M$. 

\hfill

According to \ref{ratlc}, $M$ contains a rational curve, and 
therefore, $M$ is algebraically non-hyperbolic. 
This finishes the proof of \ref{_hk_with_Lagra_non_h_Proposition_}.
\end{proof}

\hfill

\hfill

\lemma \label{_infinite_orbit_Lemma_}
Let $M$ be a projective hyperk\"ahler manifold with infinite automorphism 
group $\Gamma$. Consider the natural map $f: \Gamma \arrow \Aut(H^{1,1}(M))$. 
Then the elements of the K\"ahler cone have infinite orbits with respect to 
$f(\Gamma)$. 

\hfill

\begin{proof}
See the discussion in section 2 of \cite{_O:bimer_}. 
\end{proof}

\hfill

\lemma \label{_fin_im_Lemma_}
Let $M$ be a projective hyperk\"ahler manifold, and $\Gamma$ its 
automorphism  group. Consider the natural map 
$g: \Gamma \arrow \Aut(H_{tr}^2(M)) \times \Aut(H^{1,1}(M))$. Then $g(\Gamma)$ 
is finite in the first component $\Aut(H_{tr}^2(M))$. 

\hfill

\begin{proof}
This has been proven by Oguiso, see \cite{_Oguiso:Aut_}. 
The idea is that the BBF form restricted to the transcendental part 
$H_{tr}^2(M)$ is of K3-type. Then we can apply Zarhin's theorem 
(Theorem 1.1.1 in \cite{_Zar_}) to deduce that $g(\Gamma) \subset 
\Aut(H_{tr}^2(M))$ is finite. 
\end{proof}

\hfill

\theorem \label{Aut}
Let $M$ be a projective hyperk\"ahler manifold with infinite automorphism 
group. Then $M$ is algebraically non-hyperbolic. 

\hfill

\begin{proof} 
In the notations introduced above, 
for any K\"ahler class $w$ on $M$, its $f(\Gamma)$-orbit is infinite by 
\ref{_infinite_orbit_Lemma_}.
Fix a polarization $w$ on $M$ with normalization $q(w)=1$. 
For a given constant $C>0$ consider the set 
$${\cal D}_C = \{ x\in H^{1,1}(M, \Z)\ \  |\ \  q(x) \geq 0,\ \  
q(x,w) \leq C \}.$$
Notice that ${\cal D}_C$ is compact. 
Indeed, $y= x - q(x,w)w$ is orthogonal to $w$ with respect to the BBF 
form $q$. The quadratic form $q$ is of type $(1, \rho - 1)$ on 
$H^{1,1}(M, \Z)$ and since $q(w)>0$, the restriction 
$q|_{w^\perp}$ is negative-definite. A direct computation shows that 
$q(y)=q(x) - 2q(x,w)^2 + q(x,w)^2 q(w) = q(x) -q(x,w)^2 \geq -C^2$. 
The set ${\cal D}_C$ is equivalent to the set of elements 
$\{y \in w^\perp| q(y) \geq -C^2 \}$, which is compact because 
$q|_{w^\perp}$ is negative-definite. 
Since the set ${\cal D}_C$ is compact, 
$\sup_{x \in \Gamma \cdot \eta} \deg x = \infty$, which means 
there is a class of a curve $\eta$ with $q(\eta)>0$. However, all 
curves in the orbit $\Gamma \cdot \eta$ have constant genus. 
Since their degrees could be arbitrarily high, then $M$ is algebraically 
non-hyperbolic. 
\end{proof}

\hfill

\lemma \label{rat_curve}
Let $M$ be a hyperk\"ahler manifold such that the positive cone 
does not coincide with the K\"ahler cone. Then $M$ contains a 
rational curve. 

\hfill

{\bf Proof:} 
This is a classical result that Boucksom and Huybrechts knew in the 
early 2000's 
\cite{_Boucksom-cone_,_Huybrechts:cone_}. \endproof

\hfill

\theorem
Let $M$ be a hyperk\"ahler manifold with Picard rank $\rho$. 
Assume that either $\rho > 2$ or $\rho =2$ and the SYZ conjecture holds. 
Then $M$ is algebraically non-hyperbolic. 

\hfill

\begin{proof}
Notice that the Hodge lattice $H^{1,1}(M,\Z)$ of a hyperk\"ahler manifold
has signature $(1,k)$. Therefore, for $\rho \geq 2$, the Hodge lattice
contains a vector with positive square, and $M$ is projective
(\cite{Huy}). 
First, consider the case when $\rho > 2$. If the K\"ahler cone 
coincides with the positive cone, then the automorphism group 
$\Aut(M)$ is commensurable with the group of 
isometries $SO(H^2(M, \Z))$ (Theorem 2.17 in \cite{_AV:Aut_})
preserving the Hodge type. By \ref{_fin_im_Lemma_}, this group
is commensurable with the group of isometries of the Hodge lattice
$H^{1,1}(M,\Z)$.
By Borel and Harish-Chandra's theorem (\cite{BHC}), if $\rho>2$,
any arithmetic subgroup of $SO(1, \rho -1)$ is a lattice.
However, Borel density theorem implies that any lattice in a non-compact simple Lie group 
is Zariski dense (\cite{_Borel:density_}). 
Therefore, for $\rho>2$, $SO(H^{1,1}(M, \Z))$ is infinite. 
In this case $\Aut(M)$ is also infinite and 
we can apply \ref{Aut}. On the other hand, if the K\"ahler cone does not 
coincide with the positive cone, then by \ref{rat_curve} there 
is a rational curve on $M$. Therefore, $M$ is algebraically non-hyperbolic. 

Now let $\rho =2$. Assume the positive cone and the K\"ahler cone 
coincide. If there is no $\eta \in H^{1,1}(M, \Z)$ with $q(\eta)=0$, 
then by Theorem 87 in 
\cite{Dickson}, $SO(H^{1,1}(M, \Z))$ is isomorphic to $\Z \times \Z/2\Z$. 
Therefore, both $SO(H^{1,1}(M, \Z))$ and $\Aut(M)$ are infinite and 
we can apply \ref{Aut}. 
If there is $\eta \in H^{1,1}(M, \Z)$ with $q(\eta)=0$, then 
the SYZ conjecture implies that $\eta$ defines a rational fibration 
on $M$ and we could apply \ref{_hk_with_Lagra_non_h_Proposition_}. 
If $\rho=2$ and the positive and the K\"ahler cones are 
different (i.e., the positive cone is divided into K\"ahler chambers), 
then there is a nef class $\eta \in H^{1,1}(M, \Z)$ with 
$q(\eta)=0$. Since we assumed that the SYZ conjecture holds, the class 
$\eta$ defines a Lagrangian fibration on $M$. Applying 
\ref{_hk_with_Lagra_non_h_Proposition_} we conclude that $M$ is 
algebraically non-hyperbolic. 
\end{proof}

\hfill

\remark We conjecture that all projective hyperk\"ahler manifolds
are algebraically non-hyperbolic. However, our proof fails for manifolds
with Picard rank 1.

\hfill

{\bf Acknowledgments.} This work was inspired by a question of Erwan 
Rousseau about algebraic non-hyperbolicity of hyperk\"ahler manifolds. 
The paper was started while the first-named author was visiting 
Universit\'e libre de Bruxelles and she is grateful to Joel Fine for 
his hospitality. It was finished at the SCGP during the second-named author's 
stay there. We are grateful to the SCGP for making this possible.

\noindent {\sc Ljudmila Kamenova\\
Department of Mathematics, 3-115 \\
Stony Brook University \\
Stony Brook, NY 11794-3651, USA,} \\
\tt kamenova@math.sunysb.edu
\\

\noindent {\sc Misha Verbitsky\\
{\sc Laboratory of Algebraic Geometry,\\
National Research University HSE,\\
Faculty of Mathematics, 7 Vavilova Str., \\Moscow, Russian Federation,}\\
\tt  verbit@mccme.ru}, also: \\
{\sc Universit\'e libre de Bruxelles, CP 218,\\
Bd du Triomphe, 1050 Brussels, Belgium
}

\end{document}